\documentclass[12pt]{article}
\begin{document}
\parindent=0pt
\centerline{\Large{Painlev\'e's theorem extended}}
\vskip0,2truecm
\centerline{\large{Claudio Meneghini}}
\vskip0,3 truecm
{
\small
{\bf Abstract:}
we extend 
{ Painlev\'e's determinateness theorem }
from the theory of ordinary differential equations in the complex domain
allowing more general 
'multiple-valued'
Cauchy's problems.
We study $C^0-$continuability 
(near singularities) of solutions.
} 
\bibliographystyle{plain} 

\def\loft{}
\def\roght{}
\def\QUAN{\vrule height6pt width6pt depth0pt}
\def\blof #1{\vskip0,2truecm{\bf #1}\vskip0,2truecm }
\newtheorem{definition}{Definition}
\newtheorem{lemma}[definition]{Lemma}
\newtheorem{proposition}[definition]{Proposition}
\newtheorem{theorem}[definition]{Theorem}        
\newtheorem{corollary}[definition]{Corollary}  
\newtheorem{remark}[definition]{Remark}

\font\sdopp=msbm10
\font\cir=wncyb10
\def\ERRE {\sdopp {\hbox{R}}}
\def\CI {\sdopp {\hbox{C}}}
\def\DI {\sdopp {\hbox{D}}}
\def\ENNE{\sdopp {\hbox{N}}}
\def\PI {\sdopp {\hbox{P}}}
\def\R{\Delta} 
\def\NORM{\hbox{\boldmath{}$\vert\vert$\unboldmath}}
\blof{Foreword and preliminaries}
In this paper we 
slightly
improve
{Painlev\'e's determinateness theorem }
(see \cite{hille}, th.3.3.1),
investigating the $C^0-$continuability 
of the solutions of
finitely {'multiple-valued'}, meromorphic
Cauchy's problems.
In particular, will shall be interested in 
phenomena taking place when the attempt of
continuating a solution along an arc leads to 
singularities of the known terms: we shall see
that, under not too restricting hypotheses, 
this process will  converge to a limit.
Of course we shall formalize 
'multiple valuedness' by means of
{Riemann domains} over regions
in $\CI^2$ (see also \cite{gunros}, p.43 ff).
Branch points will be supposed
to lie on algebraic curves.
We recall 
that in the classical statement of the theorem 
'multiple valuedness' in the
known term
is allowed with respect to the independent variable only.

The following theorem
extends to the complex domain the so called 'single-sequence criterion' from 
the theory of real o.d.e.'s (see e.g.\cite{giusti}, th. 3.2); a technical 
lemma ends the section; 
the local existence-and-uniqueness theorem is reported in the appendix.

\begin{theorem}
Let $W$ be a $\CI^N$-valued holomorphic mapping, solution of the equation $W^{\prime}(z)=F
\loft(W(z),z    
\roght)$ in ${\cal V}\subset\CI$ , where 
$F$ is a  $\CI^N$-valued holomorphic mapping
in a neighbourhood of 
$graph(W)$.
Let $z_{\infty}\in \partial{\cal V}$, suppose
that there exists a sequence $\{z_n\}\rightarrow z_{\infty}$, such that,
set $W(z_n):=W_n$,
 $\lim_{n\to\infty}W_n=W_{\infty}\in\CI^N$
and that $F$ is holomorphic at $
\loft(W_{\infty},z_{\infty}    
\roght)$:
then $W$ admits analytical continuation up to $z_{\infty}$.
\label{complseq}
\vskip0,1truecm
\rm
{\bf Proof:}
 we deal only with the case $N=1$:
we can find $a>0$ and $b>0$ such that the Taylor's developments
$\sum_{k,l=0}^{\infty}c_{kln}(W-W_n)^k(z-z_n)^l$ at $(W_n,z_n)$
of $F$ are absolutely and uniformly convergent in all closed bidiscs
$\overline{\DI((W_n,z_n),a,b)}$. 
By means of Cauchy's estimates we can find an upper bound $T$
for $\sum_{k,l=0}^{\infty}\vert c_{kln}\vert a^k b^l$; by theorem 2.5.1
of \cite{hille} the solutions $S_n$ of 
$
W^{\prime}=F
\loft(W,z    
\roght),\   
W(z_n)=W_n
$
have radius of convergence at 
least
$a(1-e^{-b/2aT}):=\sigma$.
Thus
there exists $M$ such that   $z_{\infty}\in\DI(z_M,\sigma)$; 
by continuity,
$S_{M}(z_{\infty})=W_0$,
and, by uniqueness,
$
S_{M}=
S_{\infty}
$ 
in $\DI(z_M,\sigma)\cap\DI(z_{\infty},\sigma)$,
i.e. $W$ admits analytical continuation up to $z_\infty$.
\QUAN
\end{theorem} 
\begin{lemma}
Let $X$ be a metric space, $\gamma:[a,b)\rightarrow X$ a continuous arc
and suppose that there 
does not exist $\lim_{t\to b}\gamma(t)$: then, for every 
$N-$tuple $\{x_1...x_N\}\subset X$ there exists a sequence
$\{t_i\}\rightarrow b$ and neighbourhoods
$U_k$ of $x_k$ such that $\{\gamma(t_i)\}\subset X\setminus 
\bigcup_{k=1}^N U_k$.
\label{metrico}
\QUAN
\end{lemma}
We recall that
a {\sl Riemann domain} over a region 
${\cal U}\subset\CI^N$ is a complex manifold $\R$ with
an everywhere maximum-rank holomorphic surjective mapping 
$p:\R\rightarrow{\cal U}$;
$\R$ is {\sl proper} provided that so is $p$
(see \cite{gunros} p.43).

\blof{The main theorem}
Let ${\cal N}$ 
be a curve
in $\CI^2$,
$(\R,p)$ a proper Riemann domain 
over $\CI^2\setminus{\cal N}$,
$F$ a meromorphic function on $\R$,
holomorphic outside a curve ${\cal M}$,
$X_0\in \R\setminus{\cal M}$,
$(u_0,v_0)=p(X_0)$ 
and $\eta$ 
a local inverse of $p$, defined  in a 
bidisc $\DI_1\times\DI_2 $
around $(u_0,v_0)$.
We shall 
set ${\cal A}={\cal N}\cup p({\cal M})$
and
denote by $u:\DI(v_0,r)\rightarrow
\DI_1$ (with $r$ like in the existence and uniqueness theorem
in the appendix)
the solution of
Cauchy's problem:
$
u^{\prime}(v)
=F\circ\eta(u(v),v)
$
,
$ 
u(v_0)=u_0
$.
\begin{theorem}
\label{painleve'} 
{\tt\ (Painlev\'e's determinateness theorem\ )}
Suppose that ${\cal A}$ is 
algebraic;
let
$\gamma\colon [0,1]\rightarrow\CI$
be an embedded $C^1$ arc starting at $v_0$ 
such that, for each $t\in[0,1]$, the complex line $v=\gamma(t)$
is not contained  in ${\cal A}$;
suppose that
an analytical continuation $\omega$ of $u$
may be got along 
$\gamma\vert_{[0,1)}$: then
there exists 
$
\lim_{t\to 1} \omega\circ\gamma(t)
$,
 within $\PI^1$.
\vskip0,1truecm
\rm
{\bf Proof:}
suppose, on the contrary, that such limit
does not exist: for every $\nu\in\CI_v$, set 
$
W_{\nu}=
pr_1
\loft(
{\cal A}
\cap 
(\CI\times\{\nu\})
\roght)
$: 
by 
hypothesis
$W_{v_1}$ is finite or empty.
The former case is trivial; 
as to the latter,
say, $W_{v_1}=\{\lambda_k\}_{k=1....q}$.
For each $k=1...q$ and every $\varepsilon>0$, set 
$
D_{k\varepsilon}=\DI(\lambda_k,\varepsilon)$,
$
T_{k\varepsilon}=\partial(\DI_{k\varepsilon})
$;
then there exists $\varrho_{\varepsilon}>0$
such that 
$
v\in\overline{\DI
\loft(v_1,\varrho_{\varepsilon}
\roght)}
\Rightarrow
W_{v}
\subset\bigcup_{k=1}^{q}D_{k\varepsilon}
$: set now 
$$
\cases
{
M_{\varepsilon}=\max_{X\in 
p^{-1}
\loft( \bigcup_{k=1}^{q}T_k \times 
\overline{\DI
\loft(v_1,\varrho_{\varepsilon}    
\roght)   }
\roght)
}
\vert F(X) \vert  
\cr
M_{R\varepsilon}=
\max_{X\in
p^{-1}
\loft(\partial(\DI(0,R))
\times 
\overline
{\DI
\loft(v_1,\varrho_{\varepsilon}
\roght)}
\roght)
}
\vert F(X) \vert  
& (for each $R>0$);
}
$$
Introduce the compact set 
$\Theta_{R\varepsilon}=
\overline{\DI(O,R)}\setminus\bigcup_{k=1}^{q}D_{k\varepsilon}$:
for 
every
$v\in\overline{\DI
\loft(v_1,\varrho_{\varepsilon}
\roght)}  $, $p^{-1}
\loft(\CI_u\times \{v\}    
\roght)$ is a Riemann surface, hence, by maximum principle,
and
by
the arbitrariness of $v$ in $\overline{\DI
\loft(v_1,\varrho_{\varepsilon}    
\roght)  }$,
$$
X\in p^{-1}
\loft(\Theta_{R\varepsilon}\times \overline{\DI
\loft(v_1,\varrho_{\varepsilon}    
\roght)  }
\roght)\Rightarrow \vert F(X) \vert\leq \max
\loft(M_{\varepsilon},
M_{R\varepsilon}
\roght)
.$$
Now we have assumed that $\omega\circ\gamma(t)$ does not admit limit as 
$t\rightarrow 1$, hence, by lemma
\ref{metrico} (with $X=\PI^1$,
$\{x_k\}=\{\lambda_k\}\cup\{\infty\}$),
there exist: a sequence $\{t_i\}\rightarrow 1$,
$\varepsilon$ small enough and $R$ large enough
such that $\{\omega(\{\gamma(t_i)\})\}\subset \Theta_{R\varepsilon}$.
Without loss of generality, we may suppose that $\{\gamma(t_i)\}\subset \DI
\loft(v_1,\varrho_{\varepsilon}
\roght) $.
Since $p$ is proper, $p^{-1}
\loft( \Theta_u\times\overline{\DI
\loft(v_1,\varrho   
\roght) }   
\roght)$ is compact, hence we
can
extract a convergent subsequence $\Omega_k$ from $p^{-1}\{\omega(\gamma(t_i)),\gamma(t_i)\}$,
whose limit we shall call $\Omega$.
By hypothesis there exists a holomorphic function element
$
\loft({\cal V},\widetilde \omega    
\roght)$ such that
${\cal V}\supset\gamma
\loft([0,1)    
\roght)$
and
$\widetilde \omega(v_0)=u(v_0)$;
moreover,
$
F\circ\eta$
could be analytically continuated across
$\left.\omega\circ\gamma\times\gamma\right\vert_{[0,1)}$, since $\widetilde \omega^{\prime}$
is finite at each point of $\gamma([0,1))$;
by constrution, $F$ is holomorphic at $\Omega$
and $rk
\loft(p_*(\Omega)
\roght)=2$, hence
$
F\circ\eta$
admits analytical continuation up 
to $\Omega$.
Therefore, by theorem \ref{complseq}, 
$\widetilde \omega$ admits analytical continuation up to $v_1$, hence there exists
$
\lim_{t\to 1}\omega\circ\gamma(t)=
\lim_{t\to 1}\widetilde \omega\circ\gamma(t)
$,
which is a contradiction.
\QUAN
\end{theorem} 

{\bf A simple example} 

Consider 
$
w^{\prime}
=
-\sqrt[4]{8}\sqrt{3z+w^2}
\,\hbox{\large /}\,
4\sqrt[4]{(z+w^2)^3}
$,
$
w(1)=1,
$
where the branches of the roots are those ones which
take positive values on the positive real axis
(this is in fact the choice of $\eta$); the problem is solved
by $w(z)=\sqrt{z}$, which admits analytical continuation, for example, 
on the backwards oriented semi-closed interval $[1,0)$: note that 
nor $3z+w^2=0$
nor $z+w^2=0$
contains any complex line $z=$const; 
hence, as expected, $\lim_{t\to 0^+,t\in\ERRE}\sqrt{t}$ 
exists, being in fact $0$.
\vskip0,1truecm
{\bf Appendix: the existence-and-uniqueness theorem}
\par
Let $W_0
$ be a complex $N-$tuple,
$z_0\in\CI$;
let $F$ be a $\CI^N-$valued holomorphic mapping in 
$\prod_{j=1}^N\DI
\loft(W_0^j,b    
\roght)\times\DI
\loft(z_0,a   
\roght)$,
($a,b\in\ERRE$)
with $C^0-$norm $M$ and
$C^0-$norm of each
${\partial F}/{\partial w^j}$ ($j=1..N$) not exceeding $K\in\ERRE$.

{\tt Theorem:}
if  
$r<min(a,b/M,1/K)$,
there exists a unique holomorphic mapping 
$
W\colon \DI
\loft(z_0,r    
\roght)\rightarrow
\prod_{j=1}^N\DI
\loft(W_0^j,b    
\roght)
$
such that 
$W^{\prime}=F(W(z),z)$ and
$W(z_0)=W_0$.
(see e.g. \cite{hille}, th 2.2.2, \cite{ince} p.281-284)
%

\tt
\end{document}